\documentclass[12pt]{amsart}
\usepackage{amssymb}

\newtheorem{theorem}{Theorem}[section]
\newtheorem{lemma}[theorem]{Lemma}
\newtheorem{proposition}[theorem]{Proposition}

\newtheorem{example}[theorem]{Example}
\newtheorem{definition}[theorem]{Definition}
\newtheorem{remark}[theorem]{Remark}
\newtheorem{corollary}[theorem]{Corollary}

\newcommand{\M}{\mbox{$\mathbb{M}$}}
\newcommand{\C}{\mbox{$\mathbb{C}$}}
\newcommand{\N}{\mbox{$\mathbb{N}$}}
\newcommand{\R}{\mbox{$\mathbb{R}$}}
\newcommand{\K}{\mbox{$\mathbb{K}$}}
\newcommand{\B}{\mbox{$\mathbb{B}$}}

\newcommand{\Hi}{\mbox{${\mathcal H}$}}
\newcommand{\Ml}{\mbox{${\mathcal M}$}}

\begin{document}
\title[Characterizations of essential ideals as operator modules over 
C*-algebras]{Characterizations of essential ideals as operator modules over 
C*-algebras}

\vspace{30 mm}

\author{Masayoshi Kaneda* and Vern Ival Paulsen*}
\address{Department of Mathematics\\University of Houston\\4800 Calhoun Road\\
Houston, TX 77204-3476 U.S.A. and
Graduate School of Mathematical Sciences\\University of Tokyo, 3-8-1 Komaba\\
Meguro-ku\\Tokyo 153-8914 JAPAN}

\email{kaneda@math.uh.edu}
\address{Department of Mathematics\\University of Houston\\4800 Calhoun Road\\
Houston, TX 77204-3476 U.S.A. }
\email{vern@math.uh.edu}
\date{April 25, 2000.}
\thanks{* Supported by a grant from the NSF}

\maketitle

\vspace{40 mm}

\begin{abstract}
In this paper we give characterizations of essential left ideals of a 
C*-algebra $A$ in terms of their properties as operator $A$-modules. 
Conversely, we seek C*-algebraic characterizations of those ideals $J$ in $A$ 
such that $A$ is an essential extension of $J$ in various categories of 
operator modules. In the case of two-sided ideals, we prove that all the above 
concepts coincide. We obtain results, analogous to \mbox{M. Hamana's} results, 
which characterize the injective envelope of a C*-algebra as a maximal 
essential extension of the C*-algebra, but with completely positive maps 
replaced by completely bounded module maps. By restricting to one-sided ideals,
module actions reveal clear differences which do not show up in the two-sided 
case. Throughout this paper, module actions are crucial. 
\end{abstract}


\pagebreak
\newpage

\section{Introduction.}
Let $A$ denote a C*-algebra. In the C*-algebra literature a two-sided ideal,
$J$, of $A$ is called {\em essential} if $aJ=\{ 0\} $ implies $a=0.$ Given a
category and an object $X$ contained in an object $Y,$ then $Y$ is called an
{\em essential extension} of $X$ provided that every morphism of $Y$ that
restricts to be an isomorphism of $X$ must necessarily be an isomorphism of
$Y.$ In this paper we study the relationships between these two different
approaches to essentiality as we vary the category.

Every left ideal can be regarded as an object in the category of operator left 
$A$-modules, with morphisms either the completely contractive or completely 
bounded left $A$-module maps. Thus, it is interesting to know the relationships
between an ideal $J$ being essential in $A$, in the C*-algebraic sense, and $A$
being an essential extension of $J$ in these categories. For two-sided ideals  
we show that all three concepts coincide. Moreover, these conditions are 
equivalent to the injective envelopes, $I(J)$ and $I(A)$, being completely 
isometrically isomorphic via a map that restricts to be the identity on $J.$ 

However, for left ideals these concepts differ. It is still true that if $I(J)$
and $I(A)$ are completely isometrically isomorphic via a map that restricts to
be the identity on $J,$ then $J$ is an essential left ideal in $A.$ But the 
converse does not necessarily hold.  

In this paper, we characterize the essential left ideals of $A$ via their 
properties as modules in the above categories. We also give some characterizations of the left 
ideals for which $A$ is an essential extension in some of these categories in 
C*-algebraic terms. But there are some categories for which we are unable to give C*-algebraic characterizations of these families of left ideals.

As in the earlier work of \mbox{M. Frank} and the second author \cite{FP}, we
find that when the morphisms are the completely contractive maps, then our
questions reduce to the case of completely positive maps and rely mostly on 
the earlier results of \mbox{M. Hamana}. However, when the morphisms are taken 
to be the completely bounded maps or in the case of one-sided ideals, then we 
are dealing with a completely isomorphic situation which requires new 
techniques or the recent results of \cite{FP} and \cite{BP}.

Two of \mbox{M. Hamana's} original characterizations of the injective envelope
$I(A)$ of $A$ were as the minimal injective object containing $A$ and as a 
maximal essential extension of $A$ in the category of operator systems and 
completely positive maps. In the work of \mbox{M. Frank} and the second author 
\cite{FP}, it was shown that $I(A)$ could also be characterized as the minimal 
injective containing $A$ in the category of operator left $A$-modules and 
completely bounded maps. Completing this set of ideas we show that $I(A)$ can 
also be characterized as the maximal essential extension of $A$ in this latter
category.
 
Throughout this paper, C*-algebras are not necessarily unital. In Section 
\ref{Essential extensions with module actions.} we cite a result of \cite{FP} 
where `unital' is assumed. But as is remarked there, many results, including 
this one, follow in the case of a non-unital C*-algebra $A$, by taking its 
minimal unitization $A^1$ and observing that $A$-modules are naturally 
$A^1$-modules.

\section{Definitions.}
We begin by recalling a remark from \cite{FP}. When dealing with a category 
whose morphisms are vector spaces of bounded linear maps, then requiring that 
an object is injective only guarantees that every bounded linear map has a 
bounded extension, but not necessarily of the same norm. When we want to 
insist that extensions also have the same norm, then we shall refer to the 
objects as {\it tight} injectives. It is easy to see that tight injectivity is 
equivalent, by scaling, to requiring that every contractive linear map has a 
contractive linear extension. Thus tight injectivity is the same as saying that
the object is injective in a category where the morphisms are all contractive 
linear maps. Throughout this paper we shall refer to situations as tight, or 
simply omit the word tight, instead of constantly referring to a change of 
categories. But of course, both viewpoints are useful. 

We start with 
definitions of several notions many of which are already known, but perhaps not
in our particular setting.

\begin{definition}
{\em Let $A$ and $B$ be operator algebras.}

{\em (1) An} operator $A-B$-bimodule {\em $X$ is an operator space which is
also an $A-B$-bimodule with
$\| (a_{lm})(x_{lm})\| \le \| (a_{lm})\| \| (x_{lm})\|$ and 
$\| (x_{lm})(b_{lm})\| \le \| (x_{lm})\| \| (b_{lm})\|$
$\forall (x_{lm})\in \M _n(X)$, $\forall (a_{lm})\in \M _n(A)$,
$\forall (b_{lm})\in \M _n(B)$, $\forall n\in \N $. If $A$ (resp., $B$) 
contains identity $1_A$ (resp., $1_B$), we always assume that $\forall x, 
1_Ax=x$ (resp., $\forall x, x1_B=x$) We also say an} operator 
left A-module (resp., operator right A-module, operator A-bimodule) {\em for an
operator $A-\C$-bimodule (resp., operator $\C-A$-bimodule, operator 
$A-A$-bimodule).}

{\em (2) A} (resp., tight) $A-B$-rigid extension {\em of an operator
$A-B$-bimodule $X$ is a pair $(Y,\iota)$ of an operator $A-B$-bimodule $Y$ and 
a completely bounded isomorphism\footnote{A one-to-one completely bounded 
linear map whose inverse $\iota (X)\longrightarrow X$ is also completely 
bounded.} (resp., a complete isometry) $\iota : X\longrightarrow Y$ which is 
also an $A-B$-bimodule map, such that the identity map on $Y$ is the only 
completely bounded (resp., completely contractive) $A-B$-bimodule map on $Y$ 
which fixes each element of $\iota (X).$}

{\em (3) A} (resp., tight) $A-B$-essential extension {\em of an operator 
$A-B$-bimodule $X$ is a pair $(Y,\iota)$ of an operator $A-B$-bimodule $Y$ and 
a completely boundedly isomorphic (resp., completely isometric) 
$A-B$-bimodule map $\iota : X\longrightarrow Y$ such that for any operator 
$A-B$-bimodule $W$ and for any completely bounded (resp., completely 
contractive) $A-B$-bimodule map $\phi : Y\longrightarrow W$, $\phi$ is a 
completely bounded isomorphism (resp., a complete isometry) whenever 
$\phi |_{\iota (X)}$ is. If there is no proper embedding of $Y$ by a completely
boundedly isomorphic (resp., completely isometric) $A-B$-bimodule map into any 
(resp., tight) $A-B$-essential extension of $X$ (or, equivalently, there is no 
proper (resp., tight) $A-B$-essential extension of $Y$), then we say 
$(Y,\iota)$ is a} maximal {\em (resp., tight) $A-B$-essential extension of $X$.
}

{\em (4) An operator $A-B$-bimodule $W$ is} (resp., tight) $A-B$-injective,
{\em if for any operator $A-B$-bimodules $Y$ and $Z$ with  $Y\subset Z$ ($Y$ 
has a matrix norm and a module action inherited from $Z$) any completely 
bounded $A-B$-bimodule map $Y\longrightarrow W$ extends to a completely bounded
$A-B$-bimodule map $Z\longrightarrow W$ (resp., preserving the completely 
bounded norm).}

{\em (5) A} (resp., tight) $A-B$-injective extension {\em of an operator 
$A-B$-bimodule $X$ is a pair $(Y,\iota)$ of an operator $A-B$-bimodule $Y$ and 
a completely boundedly isomorphic (resp., completely isometric) $A-B$-bimodule 
map $\iota : X\longrightarrow Y$ with $Y$ (resp., tight) $A-B$-injective.}

{\em (6) An} $A-B$-injective envelope {\em of an operator $A-B$-bimodule $X$ is
a pair $(Y,\iota)$ of an operator $A-B$-bimodule $Y$ and a complete isometric 
$A-B$-bimodule map $\iota: X\longrightarrow Y$ with $Y$ a minimal tight 
$A-B$-injective operator $A-B$-bimodule containing $\iota(X)$ (i.e. for any 
tight $A-B$-injective operator $A-B$-bimodule $Z$ with 
$\iota (X)\subset Z\subset Y$, necessarily $Y=Z$ holds, where $Z$ has a matrix 
norm and a module action inherited from $Y$.). We denote this $Y$ by $I(X)$. 
Frequently we identify $X$ with $\iota (X)$ and regard it as $X\subset I(X)$.}

{\em In the above (2), (3), (5), and (6), often we simply write $Y$ for 
$(Y,\iota).$}

\end{definition}

\begin{remark}
\label{definition}
{\em (1) There are still some difficulties involved in trying to define 
injective envelopes in the non-tight setting.
In (6) above we really are defining injective envelopes only in the tight 
setting.

(2) The following is a well-known argument (\cite{H2}). Let $A$ and $B$ be 
C*-algebras. For any operator $A-B$-bimodule, take any tight $A-B$-injective
extension $(Y,\iota)$ of $X$ \footnote{Such an extension exists even if $X$ is 
an operator $A-B$-bimodule, since $X$ can be considered as $X\subset \B (\Hi )$
for some Hilbert space $\Hi $ by the representation theorem for operator 
modules (\cite{CES}, \cite{B2}, \cite{BP}) and $\B (\Hi)$ is a tight 
$A-B$-injective by the bimodule version of \mbox{G. Wittstock's} extension 
theorem (Theorem 4.1 in \cite{Wi}. See also \cite{S}.).}. Given a minimal 
$\iota (X)$-projection on $Y$ \footnote{The existence of a minimal projection 
can be shown in the same way as \cite{H2}.}, then it follows that ``$(Y,\iota)$
is an $A-B$-injective envelope of $X$ $\Longleftrightarrow $ $Y$ is the range 
of a minimal $\iota (X)$-projection on some $A-B$-injective operator 
$A-B$-bimodule containing $\iota (X)$''. In particular, for any operator 
$A-B$-bimodule  $X$, its tight $A-B$-injective  envelope exists.

(3) When $A$ and $B$ are C*-algebras, as is proved in Corollary 2.6 in 
\cite{BP}, an operator $A-B$-bimodule is tight $A-B$-injective if and only if 
it is $\C -\C $-injective, and also when $X$ is an operator $A-B$-bimodule, $Y$
is an $A-B$-injective envelope of $X$ if and only if $Y$ is a 
$\C -\C$-injective envelope of $X$. For this reason, when $X$ is an operator 
$A-B$-bimodule we frequently call $I(X)$ just an} injective envelope {\em of 
$X$. Also this fact justifies that we use the simple notation $I(X)$ without 
mentioning $A$ and $B.$}
\end{remark}

\section{Essential extensions with module actions.}
\label{Essential extensions with module actions.}
In this section we obtain the necessary characterizations of injective 
envelopes as maximal essential extensions. In the tight situation, these are 
basically restatements of results of \mbox{M. Hamana.} However, in the 
non-tight situation, these characterizations rely on recent work of M. Frank 
and the second author (\cite{FP}) which extends Hamana's rigidity results to 
completely bounded module maps. Analogous results to the known results in the 
category of modules over a ring (Theorem 2.10.20 in \cite{R}) were obtained by 
\mbox{M. Hamana} in the cases of Banach modules, C*-algebras, operator systems 
and operator spaces (\cite{H0}, \cite{H1}, \cite{H2}, \cite{H3}). Without any 
intrinsic changes in his proofs, similar results hold for the case of operator
modules over C*-algebras. We summarize the key facts below.
\begin{theorem}
\label{tightext}
Let A, B be C*-algebras, and X an operator $A-B$-bimodule. Then the 
following are
equivalent:

\begin{itemize}
\item[(i)] $(Y,\iota)$ is a maximal tight $A-B$-essential
extension of X,
\item[(ii)] $(Y,\iota)$ is a tight $A-B$-injective and tight $A-B$-essential
  extension of X,
\item[(iii)] $(Y,\iota)$ is an $A-B$-injective envelope of X,
\item[(iv)] $(Y,\iota)$ is a tight $A-B$-injective and tight $A-B$-rigid
  extension of X.
\end{itemize}
Moreover, such $Y$ is complete.
\end{theorem}

\begin{proof}
First, we suppose that $X$ is an operator $A-B$-bimodule.

That (ii)$\Longleftrightarrow$(iii)$\Longleftrightarrow$(iv) immediately
follows from the operator module over C*-algebra version of Lemma 2.11 in 
\cite{H3}, but this result relies on several earlier results. For the 
convenience of the reader, we also include an outline of the proof of this 
part.

(ii)$\implies $(iii): We only need to show minimality of $Y$. Let $Z$ be a 
tight $A-B$-injective operator $A-B$-bimodule such that 
$\iota (X)\subset Z\subset Y$. Then $id_Z$ extends to a completely contractive
$A-B$-bimodule map $\phi : Y\longrightarrow Z$. By the assumption (ii),
$\phi$ has to be a complete isometry, so that $Y=Z$.

(iii)$\implies $(iv): The same as the proof of Lemma 3.6 in \cite{H2}. When
applying this proof, note that, if $\phi$ is an $\iota (X)$-projection on $Y$
(i.e. a completely contractive $A-B$-bimodule map on $Y$ with
$\phi ^2=\phi$ and $\phi |_{\iota (X)}=id_{\iota (X)}$), then
Im$\phi \subset Y$ is a tight $A-B$-injective operator $A-B$-bimodule, so that
Im$\phi$=$Y$ by minimality of $Y$. Hence $\phi (=id_Y)$ is a minimal
$\iota (X)$-projection on $Y$ and $Y$ is its range.

(iv)$\implies $(ii): The same as the proof (Necessity) of Theorem 3.7. in
\cite{H2}.

(i)$\implies$(ii): Let $(I(Y),\kappa)$ be an $A-B$-injective envelope of
$Y$. By (iii)$\implies$(ii), $(I(Y),\kappa)$ is a tight $A-B$-essential
extension of $Y$, so that $(I(Y),\kappa \circ \iota)$ is a tight
$A-B$-essential extension of $X$, hence $I(Y)=\kappa (Y)$ by maximality. Thus
$Y$ is tight $A-B$-injective.

(ii)$\implies$(i): Let $(Z,\kappa)$ be a tight $A-B$-essential extension of
$Y$. By injectivity of $Y$, $\kappa ^{-1} : \kappa (Y)\longrightarrow Y$
extends to a completely contractive $A-B$-bimodule map
$\phi : Z\longrightarrow Y$. But by essentiality of $Z$, $\phi $ has to be a
complete isometry. Thus $\kappa (Y)=Z$.

Finally, we show completeness. Let $Y$ satisfy the equivalent conditions
(i)---(iv) and let $(\widetilde{Y},\kappa)$ be its completion. Then
$\widetilde{Y} $ is an $A-B$-bimodule and obviously a tight $A-B$-essential
extension of $Y$. By injectivity of $Y$, $\kappa (Y)\longrightarrow Y$ extends
to a completely contractive $A-B$-bimodule map
$\phi : \widetilde{Y} \longrightarrow Y$ with $\phi |_{\kappa (Y)}$ completely
isometric. But by an essentiality of $\widetilde{Y} $, $\phi $ has to be a
complete isometry, so that $\widetilde{Y} =\kappa (Y)$.
\end{proof}

\begin{remark}
\label{injective envelopes}
{\em (1) From this theorem, in particular, for any operator $A-B$-bimodule
$X$, a maximal tight $A-B$-essential extension exists, since $A-B$-injective 
envelope exists (Remark\ref{definition} (2)).

(2) It immediately follows from rigidity in the above theorem that the
$A-B$-injective envelope is unique up to completely isometric
$A-B$-bimodule isomorphisms.

(3) In particular, if $X=A$ with $A$ a unital C*-algebra, $I(A)$ is completely
isometrically isomorphic to \mbox{M. Hamana's} injective envelope C*-algebra
$I_{C^*}(A)$ (\cite{H1}) of $A$. This is easily seen by observing unital
completely positive maps are just the same as unital completely contractive
maps, and using rigidities of $I(A)$ and $I_{C^*}(A)$. If $A$ is a non-unital 
C*-algebra, still $I(A)$ is completely isometrically isomorphic to 
$I_{C^*}(A^1),$ via a map that restricts to be the identity on $A,$ by 
Proposition 2.8 in \cite{BP}, where $A^1$ is the minimal unitization of $A$. 
Namely, $I(A^1)$ is completely isometrically isomorphic to $I(A)$. Actually the
above identification of $I(A)$ with $I_{C^*}(A)$ or $I_{C^*}(A^1)$ is also as 
$A-A$-, $A-\C $- and  $\C -A$-bimodules by Remark\ref{definition} (3).}
\end{remark}
 
In the non-tight case a similar result holds
thanks to Corollary 2.2 in \cite{FP} which is the generalization of `rigidity'
to the non-tight case in the presence of a module action.

\begin{theorem}
\label{ext}
Let A be a C*-algebra, I(A) its injective envelope and let E be an operator 
left $A$-module with A $\subset $ E. Then the following are equivalent:
\begin{itemize}
\item[(i)] E is a maximal $A-\C$-essential extension of A,
\item[(ii)] E is an $A-\C$-injective and $A-\C$-essential extension of A,
\item[(iii)] E is a minimal $A-\C$-injective extension of A,
\item[(iv)] E is completely boundedly isomorphic to I(A) as left $A$-modules, 
\\ via a map that restricts to the identity on A, 
\item[(v)] E is an $A-\C$-injective and $A-\C$-rigid extension of A.
\end{itemize}
The analogous result holds for the operator right A-module case and the 
operator A-bimodule case.
\end{theorem}

\begin{proof} Throughout the proof, without loss of generality, we may assume
that $I(A)=I_{C^*}(A)$ (Remark\ref{injective envelopes} (3)).

That (iii)$\implies$(iv) was proven in Theorem 3.3 in \cite{FP}.

That (iv)$\implies$(iii) follows by noting that any $A-\C$-injective
extension of $A$ which is completely boundedly isomorphic to a minimal
$A-\C$-injective extension as left $A$-modules is also minimal and that $I(A)$
is a minimal $A-\C$-injective extension by Theorem 3.1 in \cite{FP}.

(iv)$\implies$(i): Without loss of generality we may assume that
$E=I(A)$. Let $W$ be any operator left $A$-module and
$\phi : I(A)\longrightarrow W$ a completely bounded left $A$-module map which
is a completely bounded isomorphism on $A$. Let
$\psi : \phi (A)\longrightarrow A$ be the inverse of this map and extend
$\psi $ to a completely bounded left $A$-module map
$\widetilde{\psi } : W\longrightarrow I(A)$. Since
$\widetilde{\psi } \circ \phi : I(A)\longrightarrow I(A)$ is the identity on
$A$, it is the identity on $I(A)$ by rigidity (Corollary 2.2 in \cite{FP}).
This shows that $\phi (I(A))$ is completely boundedly isomorphic to $I(A)$.
Maximality: Suppose that $(F,\iota)$ is an $A-\C $-essential extension of
$I(A)$. A completely boundedly isomorphic left $A$-module map
$\iota ^{-1} : \iota (I(A))\longrightarrow I(A)$ extends to a completely
bounded left $A$-module map $\rho : F\longrightarrow I(A)$. By essentiality of 
$F$, $\rho $ has to be a completely bounded isomorphism, so that 
$F=\iota (I(A))$.

(i)$\implies$(iv): We may assume that $A\subset E\subset \B (\Hi)$ for some 
Hilbert space $\Hi$ by the representation theorem for operator modules
(\cite{CES}, \cite{B2}, \cite{BP}). Extend the identity map on $A$ to a
completely bounded left $A$-module map
$\psi :   \B(\mathcal H)\longrightarrow I(A)$, then $\psi (E)\subset I(A)$
is completely boundedly isomorphic to $E$. The map
$\psi ^{-1} : \psi (E)\longrightarrow E$ extends to a completely bounded left
$A$-module map $\phi : I(A)\longrightarrow   \B(\Hi).$ By rigidity 
(Corollary 2.2 in \cite{FP}), $\psi \circ \phi$ is the identity on $I(A)$. 
Hence we have $A\subset E\subset \phi (I(A))$ with $\phi (I(A))$ completely 
boundedly isomorphic to $I(A)$. This makes $\phi (I(A))$ an operator 
$A-\C$-essential extension of $A$ by applying (iv)$\implies$(i) that we have 
already proved. Thus, by maximality of $E$, $E=\phi (I(A))$ so that $I(A)$ is 
completely boundedly isomorphic to $E$ as left $A$-modules by $\phi$.

Now we have that (i)$\implies$(ii), since clearly (i) and (iv) imply (ii).

Next, we show that (ii)$\implies$(iv). By the assumption (ii), the identity
map on $A$ extends to a completely bounded left $A$-module map
$\phi : E\longrightarrow I(A)$ and the map $\phi (E)\longrightarrow E$ extends
to a completely bounded left $A$-module map $\psi : I(A)\longrightarrow E$.
Then $\psi \circ \phi$ is the identity on $E$. But by rigidity (Corollary 2.2
in \cite{FP}), $\phi \circ \psi$ is the identity on $I(A)$.

(iv)$\implies$(v): Without loss of generality, we may assume that $E=I(A)$.
Then (v) immediately follows from Corollary 2.2 in \cite{FP}.

(v)$\implies$(iv): By injectivity of $E$, the identity map on $A$ extends to a
completely bounded left $A$-module map $\phi : I(A)\longrightarrow E$.
Similarly, by injectivity of $I(A)$, the identity map on $A$ extends to a
completely contractive $A-B$-bimodule map $\psi : E\longrightarrow I(A)$. But
$\psi \circ \phi =id_{I(A)}$ and  $\phi \circ \psi =id_E$ by rigidity of 
$I(A)$
(Corollary 2.2 in \cite{FP}) and $E$, respectively.
\end{proof}

\vspace{5mm}

As is seen in the above proof, Corollary 2.2 in \cite{FP} played a crucial
role.

Before moving on to the main results, let us introduce some definitions.

\section{Essential left ideals.}
\label{Essential left ideals.}
Recall that a closed two-sided ideal $K$ in a C*-algebra $A$ is $essential$,
if $aK=\{ 0\}$ implies $a=0$, or, equivalently, $Ka=\{ 0\}$ implies $a=0$, or,
equivalently, $K\cap K'\neq \{ 0\}$ for all non-zero closed two-sided ideal 
$K'$ in $A$, where $a\in A$. We now generalize these ideas to the one-sided 
case. Note that a subset $J\subset A$ is a closed left ideal in $A$ if and only
if $J^*$ is a closed right ideal in $A$. Let $J$ be a closed left ideal in $A.$
As is well known (say, \cite{KR}, \cite{M}, \cite{Ta}), $J$ has a contractive
right approximate identity $\{ e_{\alpha }\}$. Namely,
$\{ e_{\alpha }\} \subset (J\cap J^*)_+$ is an increasing net and
$\lim_{\alpha \to \infty }je_{\alpha }=j \; \forall j\in J$ and
$\| e_{\alpha }\| \le 1 \; \forall \alpha$. Clearly, it is also a contractive
left approximate identity of $J^*$.

\begin{lemma}
\label{rl}
Let A be a C*-algebra, I a closed right ideal in A, and J a closed left ideal 
in A. Then $$I\cap J=IJ,$$ where 
$IJ:=\overline{\text{span}}\{ij; i\in I, j\in J\},$ here the closure is taken 
in $A.$
\end{lemma}

\begin{proof}
$I\cap J\supset IJ$ is clear. We show $I\cap J\subset IJ$. Let $i\in I\cap J$
and $\{e_{\alpha }\} $ a right approximate identity of $J$, then
$i=\lim_{\alpha \to \infty} ie_{\alpha }\in IJ$.
\end{proof}

\begin{proposition}
\label{ess}
Let A be a C*-algebra and J a closed left ideal in A. Then the following are
equivalent.
\begin{itemize}
\item[(i)] JJ* is an essential two-sided ideal in A,
\item[(ii)] If $a\in A$ and $aJ=\{ 0\} $, then a=0,
\item[(iii)] $I\cap J\neq \{ 0\} $ for any non-zero closed right ideal I in A,
\item[(iv)] $K\cap J\neq \{ 0\} $ 
for any non-zero closed two-sided ideal K in A.
\end{itemize}
\end{proposition}

\begin{proof}
(i)$\implies $(ii): Suppose $aJ=\{ 0\}$ for $a\in A$, then $aJJ^*=\{ 0\} $, so
that $a=0$.

(ii)$\implies $(iii): Let $I$ be a closed right ideal in $A$ with
$I\cap J=\{ 0\} $, then $IJ=\{ 0\} $ by \mbox{Lemma\ref{rl},} so that
$I=\{ 0\} $ by the assumption (ii).

(iii)$\implies $(iv): Clear.

(iv)$\implies $(i): Obviously $JJ^*$ is a closed two-sided ideal in $A$, so we
only show essentiality. Let $K:=\{ a\in A ; aJJ^*=\{ 0\} \} $, then $K$ is a
closed two-sided ideal in $A$. Suppose that $a\in K$ and $aJJ^*=\{ 0\} $, then
$0=ajj^*a^*=aj(aj)^*$ for all $j\in J$, so that $aJ=\{ 0\} $, hence
$K\cap J=KJ=\{ 0\} $. The assumption (iv) yields $K=\{ 0\} $.
\end{proof}

\begin{remark}
\label{remark ideals}
{\em Let $A$ and $J$ be the same as the assumptions of Proposition\ref{ess}, 
and let $\{ e_{\alpha }\}$ be a right approximate identity of $J.$ 

(1) $JJ^*$ is the two-sided ideal in $A$ generated by $J$ in $A$. Also
$JJ^*=JA$ since, for any $a\in A$,
$ja=\lim_{\alpha \to \infty } (je_{\alpha})a=j\cdot \lim_{\alpha \to \infty }
e_{\alpha }a$.

(2) Any $j\in J$ is written as $j=\lim_{\alpha \to \infty} je_{\alpha }$, hence
$J\subset JJ^*$ even if $A$ is non-unital. Similarly, $J^*\subset JJ^*$.}
\end{remark}

\begin{definition}
\label{essentialleftideals}
{\em Let $A$ and $J$ be the same as the assumptions in Proposition\ref{ess}.
When the equivalent statements in Proposition\ref{ess} hold, we say $J$ is an}
essential left ideal {\em in $A$. } Essential right ideals {\em are defined 
similarly.}
\end{definition}

\section{Main results.}
\label{Main results.}
In this section we present our main results, showing the relationships between 
the various notions of ``essential'' for left ideals. We apologize to the 
reader in advance for the rather long list of conditions that are equivalent to
a left ideal being essential in the C*-algebraic sense. However, it is 
convenient to know that all of these are equivalent. We have presented the 
results in this particular fashion so that one can more easily contrast the 
non-tight and two-sided cases. 
\begin{theorem}
\label{tight}
Let A be a C*-algebra, J a closed left ideal in A, and $\{ e_{\alpha }\} $ a
contractive right approximate identity of J. Then statements in {\em (I)} and
{\em (II)} are equivalent respectively, and ``each statement in {\em (I)}''
$\implies $ {\em (4)} $\implies $ {\em (5)} $\implies $ ``each statement in
{\em (II)}''.

\begin{tabular}{rp{6in}}
{\em (I)} &
$\begin{cases}
\text{\em (1)} & \text{I(J) is completely isometrically isomorphic to I(A),} \\
 & \text{via a map that restricts to the identity on J,} \\
\text{\em (2)} & \text{A is a tight $\C- \C$-essential extension of J,} \\
\text{\em (3)} & \text{A is a tight $A-\C$-essential extension of J,}
\end{cases} $ \\
& $\; \; \; $\text{\em (4)} $\; \; \; \;
$\text{$\| (a_{lm})\| =\sup_{\alpha }\| (a_{lm}e_{\alpha })\| $
for all $(a_{lm})\in \M _n(A)$ and for all $n\in \N $,} \\
& $\; \; \; $\text{\em (5)} $\; \; \; \; $\text{There exists $n\in \N $ such
that} \\
& $\; \; \; \; \; \; \; \; \; \; \; \; \; \| (a_{lm})\| =\sup_{\alpha }\|
(a_{lm}e_{\alpha })\| $
\text{for all $(a_{lm})\in \M _n(A)$,} \\
{\em (II)} &
$\begin{cases}
\text{\em (6)} &
\text{$\| (a_{lm})\| =\sup \{ \| (a_{lm})(j_{lm})\| ;
(j_{lm})\in \M _n(J), \| (j_{lm})\| \le 1\} $} \\
& \text{for all $(a_{lm})\in \M  _n(A)$ and for all $n\in \N $,} \\
\text{\em (7)} &  \text{There exists $n\in \N$ such that} \\
& \text{$\| (a_{lm})\| =\sup \{ \| (a_{lm})(j_{lm})\| ;
(j_{lm})\in \M _n(J), \| (j_{lm})\| \le 1\} $} \\
& \text{for all $(a_{lm})\in \M _n(A)$,} \\
\text{\em (8)} & \text{There exists $c\ge 1$ such that} \\
& \text{$\| (a_{lm})\|
\le c\cdot \sup \{ \| (a_{lm})(j_{lm})\| ;
(j_{lm})\in \M _n(J), \| (j_{lm})\| \le 1 \} $} \\
& \text{for all $(a_{lm})\in \M  _n(A)$ and for all $n\in \N $,} \\
\text{\em (9)} & \text{There exist $c\ge 1$ and $n\in \N $ such that} \\
& \text{$\| (a_{lm})\|
\le c\cdot \sup \{ \| (a_{lm})(j_{lm})\| ;
(j_{lm})\in \M _n(J), \| (j_{lm})\| \le 1 \} $} \\
& \text{for all $(a_{lm})\in \M _n(A)$,} \\
\text{\em (10)} & \text{A is a tight $\C -J$-essential extension of J,} \\
\text{\em (11)} & \text{A is a $\C -J$-essential extension of J,} \\
\text{\em (12)} & \text{A is a tight $A-J$-essential extension of J,} \\
\text{\em (13)} & \text{A is an $A-J$-essential extension of J,} \\
\text{\em (14)} & \text{$\M _{n}(J)$ is an essential left ideal in
$\M _{n}(A)$ for all $n\in \N$,} \\
\text{\em (15)} & \text{$\M _{n}(J)$ is an essential left ideal in
$\M _{n}(A)$ for some $n\in \N $,} \\
\text{\em (16)} & \text{$\M _n(A)$ is canonically *-isomorphically embedded in
$\Ml (\M _n(JJ^*))$} \\
& \text{for all $n\in \N $,} \\
\text{\em (17)} & \text{$\M _n(A)$ is canonically *-isomorphically embedded in
$\Ml (\M _n(JJ^*))$} \\
& \text{for some $n\in \N $,}
\end{cases} $
\end{tabular}
where $\Ml (\M _n(JJ^*))$ is the multiplier algebra of $\M _n(JJ^*).$
\end{theorem}

\begin{proof}
(1)$\implies$(2): Without loss of generality, we may assume that
$A\subset I(J)=I(A)$. Let $W$ be any operator space, 
$\phi : A\longrightarrow W$ any complete contraction which is completely 
isometric on $J$. Since $I(J)$ is tight $\C -\C $-injective, 
$\phi ^{-1}: \phi (J)\longrightarrow J$ extends to a complete contraction 
$\psi : W\longrightarrow I(J)$. Then, again, by tight 
$\C -\C $-injectivity of $I(J)$,  $\psi \circ \phi$ extends to a complete
contraction $\rho : I(J)\longrightarrow I(J)$ with $\rho |_J=id_J$. By rigidity
(Theorem\ref{tightext}), $\rho =id_{I(J)}$, hence $\phi $ has to be a complete 
isometry.

(2)$\implies$(3): Clear.

(3)$\implies$(1): Since $I(J)$ is tight $A-\C $-injective, $id_J$ extends to a
completely contractive left $A$-module map $\phi : I(A)\longrightarrow I(J)$.
By the assumption (2), $\phi |_A$ is a complete isometry. Since $I(A)$ is tight
$A-\C $-injective, $(\phi |_A)^{-1}: \phi (A)\longrightarrow A$ extends to a
completely contractive left $A$-module map $\psi : I(J)\longrightarrow I(A)$.
Then $\psi \circ \phi $ and $\phi \circ \psi$ are completely contractive left
$A$-module map, so that $\psi \circ \phi =id_{I(A)}$ and
$\phi \circ \psi =id_{I(J)}$ by rigidity (Theorem\ref{tightext}). Hence $\phi $
and $\psi $ are onto complete isometries, so that $I(J)$ is completely 
isometrically isomorphic to $I(A).$

(2)$\implies$(4): Let us define $\wp : A\longrightarrow \R _+$ by
$\wp (a):=\sup _{\alpha }\| ae_{\alpha }\|,$ $\forall a\in A$. And let
$\widetilde A :=A/$Ker$\wp $, then $\widetilde A$ is an operator space with a
well-defined matrix norm $\|| \cdot \||,$ where
$\|| (a_{lm}+$Ker$\wp)\|| :=\sup _{\alpha }\| (a_{lm}e_{\alpha })\|,$
$\forall (a_{lm})\in \M _n(A)$, $\forall n\in N$. Define
$\phi : A\longrightarrow \widetilde A$ by $\phi (a):=a+$Ker$\wp $. Then
$\phi $ is a completely contractive liner map which is completely isometric on 
$J$, so that, $\phi $ is completely isometric on $A$ by the assumption (2). 
Hence 
$\| (a_{lm})\| =\|| (a_{lm})\|| =\sup _{\alpha }\| (a_{lm}e_{\alpha })\|,$
$\forall (a_{lm})\in \M _n(A)$, $\forall n\in N$.

That (4)$\implies $(5)$\implies $(7)$\implies $(9)$\implies $(15), that 
(6)$\implies $(7), and that (6)$\implies $(8)$\implies $(9) are clear.

(6)$\implies $(10): Let $W$ be any operator right $J$-module and
$\phi : A\longrightarrow W$ any completely contractive right $J$-module map
which is completely isometric on $J$. Then $\forall n\in \N$ and
$\forall (a_{lm})\in \M _n(A)$,
\begin{center}
\begin{tabular}{rl}
& $\| (a_{lm})\| $ \\
= & $\sup \{ \| (a_{lm})(j_{lm})\| ;
(j_{lm})\in \M _n(J), \| (j_{lm})\| \le 1\} $ \\
= & $\sup \{\| \phi _n((a_{lm})(j_{lm}))\| ;
(j_{lm})\in \M _n(J), \| (j_{lm})\| \le 1\} $ \\
= & $\sup \{ \| (\phi (a_{lm}))(j_{lm})\| ;
(j_{lm})\in \M _n(J), \| (j_{lm})\| \le 1\} $\\
$\le $ & $\| (\phi (a_{lm}))\| $,
\end{tabular}
\end{center}
which shows that $\phi $ is completely isometric on $A.$ Where 
$\phi _n: \M _n(A)\longrightarrow \M _n(W)$ is defined by 
$\phi _n((a_{lm})):=(\phi(a_{lm})).$

(10)$\implies $(12): Clear.

(12)$\implies $(15): Let $K:=\{ a\in A; aJ=\{ 0\} \}$, then $K$ is a closed
two-sided ideal in $A$ and $K\cap JJ^*=KJJ^*=\{ 0\} $. Let $A/K$ be equipped
with the quotient norm that makes $A/K$ a C*-algebra, hence an operator
$A-J$-bimodule with induced matrix norms and natural module actions. Let
$\pi : A\longrightarrow A/K$ be the quotient map that is a *-homomorphism and
also a completely contractive $A-J$-bimodule map. Since $K\cap JJ^*=\{ 0\} $,
$\pi $ is one-to-one on $JJ^*$, hence *-isometric on $JJ^*$. Thus $\pi $ is
completely isometric on $JJ^*$, especially on $J$. By the assumption (12),
$\pi $ is one-to-one on $A$, so that $K=\{ 0\} $, which means $J$ is an
essential left ideal.
 
(6)$\implies $(11): Let $W$ be any operator right $J$-module and
$\phi : A\longrightarrow W$ any completely bounded right $J$-module map which
is completely boundedly isomorphic on $J$. Then there exists $c>0$ such that
$\| (j_{lm})\| \le c\cdot \| (\phi (j_{lm}))\|,$ 
$\forall (j_{lm})\in \M _n(J)$ with $\| (j_{lm})\| \le 1$, $\forall n\in \N $.
Hence $\forall n\in \N$ and $\forall (a_{lm})\in \M (A)$,
\begin{center}
\begin{tabular}{rl}
& $\| (a_{lm})\| $ \\
= & $\sup \| (a_{lm})(j_{lm})\| ;
(j_{lm})\in \M _n(J), \| (j_{lm})\| \le 1\} $ \\
$\le $ & $c\cdot \sup \{\| \phi _n((a_{lm})(j_{lm}))\| ;
(j_{lm})\in \M _n(J), \| (j_{lm})\| \le 1\} $ \\
= & $c\cdot \sup \{ \| (\phi (a_{lm}))(j_{lm})\| ;
(j_{lm})\in \M _n(J), \| (j_{lm})\| \le 1\} $ \\
$\le $ & $c\cdot \| (\phi (a_{lm}))\|$.
\end{tabular}
\end{center}
Thus, $\phi $ is a completely bounded isomorphism.

(11)$\implies $(13): Clear.

(13)$\implies $(15): Just the same as that (12)$\implies$(15).

(15)$\implies $(14): Suppose $aJ=\{ 0\},$ $\forall a\in A$. Let
$(a_{lm})\in \M _n(A)$ be the matrix such that $a_{11}=a$ and all other entries
are $0$. Then, by the assumption, $(a_{lm})\cdot \M _n(J)=\{ 0\} $ which 
implies $(a_{lm})=(0)$ from (15), so that $a=0$. Hence $J$ is an essential left
ideal in $A,$ so that $\M _{n'}(J)$ is an essential left ideal in 
$\M _{n'}(A),$ $\forall n'\in \N.$

(14)$\implies $(16): Fix $n \in \N.$ Since $\M _n(J)$ is an essential left
ideal in $\M _n(A)$, the canonical embedding
$\varphi : \M _n(A)\longrightarrow \Ml (\M _n(JJ^*))$ is one-to-one
(Proposition\ref{ess}), hence it is a *-isomorphism.

(16)$\implies $(17): Clear.

(16)$\implies $(6): Let an $n\in \N.$ Let us consider the following maps.
$$\M _n(A)\stackrel{\varphi}{\hookrightarrow} \Ml (\M _n(JJ^*))
\stackrel{\tilde{\rho}}{\stackrel{\sim}{\rightarrow}} \Ml (\K (\M _n(J)))
\stackrel{\tau}{\stackrel{\sim}{\rightarrow}} \B (\M _n(J))$$

Where we regard $\M _n(J)$ as a right Hilbert C*-module over $\M _n(J^*J)$ with
the inner product 
$<\cdot |\cdot > : \M _n(J)\times \M _n(J) \longrightarrow \M _n(J^*J)$ defined
by $<j|j'>:= j^*j',$ $\forall j, \forall j'\in  \M _n(J).$ 

$\B (\M _n(J))$ is the set of adjointable maps on $\M _n(J).$

$\K (\M _n(J))$ is the set of compact adjointable operators on $\M _n(J).$ 
Namely, $\K (\M _n(J)):=\overline{\text{span}}\{ \theta _{j,j'}; 
j,j'\in \M _n(J) \},$ 
where $\theta _{j,j'} : \M _n(J)\longrightarrow \M _n(J)$ is defined by 
$\theta _{j,j'}(j'') := j{j'}^*j''$ 
for $\forall j, \forall j', \forall j'' \in \M _n(J),$ and the closure is taken
in $\B (\M _n(J)).$  

The map $\varphi$ is the canonical embedding.

The map $\rho : \M _n(JJ^*)\longrightarrow \K (\M _n(J))$ is defined in the 
following way:
$\rho(\sum_m (j_m{j'_m}^*))=\sum_m \theta_{j_m,j'_m},$ 
$\forall j_m, \forall j'_m \in \M _n(J)$ is a well-defined
*-homomorphism and extends to an onto *-homomorphism 
$\M _n(JJ^*)\longrightarrow \K (\M _n(J)).$ By Proposition\ref{ess}, 
$\M _n(J)$ is an essential left ideal in $\M _n(JJ^*),$ thus $\rho$ is 
one-to-one and hence a *-isomorphism. Therefore, we obtain a *-isomorphism 
$\tilde{\rho} : \Ml (\M _n(JJ^*))\longrightarrow \Ml (\K (\M _n(J))).$

The map $\tau ^{-1}: \B (\M _n(J))\longrightarrow \Ml (\K (\M _n(J)))$ is 
defined in the following way:
$\tau ^{-1}(T):= (T_1,T_2),$ 
where $T_1(\theta _{j,j'}):= T\cdot \theta _{j,j'},$ 
$T_2(\theta _{j,j'}):= \theta _{j,j'} \cdot T,$ 
$\forall T\in \B (\M _n(J)).$ 
For a proof that $\tau ^{-1}$ is a *-isomorphism, see \cite{K}, \cite{W}. 

We can easily see that $\tau \circ \tilde{\rho} \circ \varphi (a)= T_a,$ where 
$T_a(j):= aj,$ $\forall a\in \M _n(A), \forall j\in \M _n(J).$ Hence (6) 
follows. 

(17)$\implies $(7): The same as that (16)$\implies $(6).
\end{proof}

\begin{remark}
{\em (1) $\sup_{\alpha }\| (a_{lm}e_{\alpha })\|$ does not depend on a
choice of contractive right approximate identities. In fact, take another
contractive right approximate identity $f_{\beta }.$ Then, by noting that
$\{ e_{\alpha }\} \subset J\cap J^*$ and $f_{\beta }$ is also a contractive
left approximate identity of $J^*$ (Remark\ref{remark ideals}),
$\lim_{\beta \to \infty }e_{\alpha }f_{\beta }=e_{\alpha }
=\lim_{\beta \to \infty }f_{\beta }e_{\alpha }, \; \forall \alpha .$ Hence,}
\begin{center}
\begin{tabular}{rl}
& $\sup_{\alpha }\| (a_{lm}e_{\alpha })\| $\\
= & $\sup_{\alpha }\lim_{\beta \to \infty }\| (a_{lm}e_{\alpha }f_{\beta })\| $
\\
= & $\sup_{\alpha }\lim_{\beta \to \infty }\| (a_{lm}f_{\beta }e_{\alpha })\| $
\\
$\le $& $\sup_{\alpha }\sup_{\beta }\|(a_{lm}f_{\beta }e_{\alpha })\| $\\
$\le $& $\sup_{\beta }\| (a_{lm}f_{\beta })\| .$
\end{tabular}
\end{center}
{\em Similarly, the other inequality holds. The similar thing holds in
Theorem\ref{non-tight}, Corollary\ref{two-sided} and Corollary\ref{JJ*}.}

{\em (2) It is interesting to note that injective envelopes are related to
tight ``$\C-\C$''-, or, ``$A-\C $''-essential extensions, while essential left 
ideals are related to tight ``$\C-J$''-, or, ``$A-J$''-essential extensions.}

{\em (3) In (4) and (5) of Theorem\ref{tight}, we can replace
$\| (a_{lm}e_{\alpha })\| $ by $\| (e_{\alpha }a_{lm})\| $. In fact,
$\| (a_{lm})\| =\| (a_{lm})^*\| =\sup_{\alpha }\| (a^*_{ml}e_{\alpha })\|
=\sup_{\alpha }\| (e_{\alpha }a_{lm})\| $. Similar considerations hold in
Theorem\ref{non-tight}, Corollary\ref{two-sided}, and Corollary\ref{JJ*}.}

{\em (4) (5) of Theorem\ref{tight} is equivalent to saying that
$\| a\| =\sup_{\alpha }\| ae_{\alpha }\| $ for all $a\in A$. Similar 
considerations hold in Theorem\ref{non-tight}, Corollary\ref{two-sided}, and
Corollary\ref{JJ*}.}

{\em (5) The statement (9) with $n=1$ is equivalent to saying that the left 
$A$-module $J$ is} c-faithful {\em in the terminology of \cite{B2}. Hence we 
saw that if we regard left ideals of $A$ as operator left $A$-modules, then $J$
being faithful\footnote{This is equivalent to saying that $J$ is an essential 
left ideal in $A$ in our definition (Definition\ref{essentialleftideals}).} and
$J$ being $c$-faithful are equivalent for any $c\ge 1.$}

{\em (6) In the proof that (16)$\implies$(6) and that (17)$\implies$(7), we 
used the Hilbert C*-module theory. We can give an alternative proof of these 
parts that looks easier but is less informative (The involution is not given 
explicitly.). First, note that, that (16)$\implies$(14) and that 
(17)$\implies$(15) easily follow. So it suffices to show that 
(14)$\implies$(6). We use an already known principle (\cite{T}, \cite{B1}) : 
Any contractive homomorphism from C*-algebra to Banach algebra is a 
*-homomorphism with a certain involution in the range. In fact, $\B (\M _n(J))$
is a Banach algebra and the canonical embedding
$\M _{n}(A)\longrightarrow \B (\M _n(J))$ is a contractive homomorphism. (14) 
implies that this embedding is one-to-one, hence a *-isomorphism, so that (6) 
follows.}
\end{remark}

For the non-tight case, we could not connect essential extensions with
injective envelopes. This is mainly because we can not say $I(J)$ is a
(non-tight) $A-\C $-rigid extension of $J$. But still a similar result holds.

\begin{theorem}
\label{non-tight}
Let A, J and $\{ e_{\alpha }\} $ be the same as the assumptions in
Theorem\ref{tight}. Then {\em (2a)}$\implies $ {\em (3a)} $\implies $
{\em (4a)} $\implies $ {\em (5a)} $\implies $ ``each statement of {\em (II)} in
Theorem\ref{tight}''.
 
\begin{itemize}
\item[(2a)] A is a $\C-\C$-essential extension of J,
\item[(3a)] A is an $A-\C$-essential extension of J,
\item[(4a)] There exists $c\ge 1$ such that \\
$\| (a_{lm})\| \le c\cdot \sup_{\alpha }\| (a_{lm}e_{\alpha })\| $
for all $(a_{lm})\in \M _n(A)$ and all $n\in \N $, 
\item[(5a)] There exist $c\ge 1$ and $n\in \N $ \\
such that $\|( a_{lm})\| \le c\cdot \sup_{\alpha }\| (a_{lm}e_{\alpha })\| $
for all $(a_{lm})\in \M _n(A)$.
\end{itemize}
\end{theorem}

\begin{proof}
(2a)$\implies $(3a): Clear.

(3a)$\implies $(4a): Similar to that (3)$\implies $(4) in Theorem\ref{tight}.

That (4a)$\implies $(5a)$\implies $``(9) in Theorem\ref{tight}'' is clear.
\end{proof}

\vspace{5mm} 

Note that the statement (2a) says just that $A=J.$ In fact, if 
$J\subsetneqq A,$ then we can take a $J$-projection $\phi$ on $A$ (i.e. a 
completely bounded linear map on $A$ with $\phi^2=\phi$ and $\phi|_J=id_J$) 
such that $J\subset$Im$\phi \subsetneqq A$ with the codimension of Im$\phi$ in 
$A$ is 1.  

Together with the right ideal versions of Theorem\ref{tight} and
Theorem\ref{non-tight}, and by noting that, when $J$ is a two-sided ideal in 
$A,$ $\C -J$- or $A-J$- in (10)---(13) of Theorem\ref{tight} can be replaced by
$\C-A$- or $A-A$- with trivial modifications in the proof, the next corollary 
immediately follows.

\begin{corollary}
\label{two-sided}
Let A be a C*-algebra, K a closed two-sided ideal in A, and $\{ e_{\alpha }\}$
a contractive approximate identity of K. Then the following are equivalent.
\begin{itemize}
\item[(1)] I(K) is completely isometrically isomorphic to I(A), \\
via a map that restricts to the identity on K,
\item[(2)] A is a tight $\C-\C$-essential extension of K,
\item[(3)] A is a tight $A-\C $-essential extension of K,
\item[(3a)] A is an $A-\C $-essential extension of K,
\item[(4)] $\| (a_{lm})\| =\sup_{\alpha }\| (a_{lm}e_{\alpha })\| $ for all
$(a_{lm})\in \M _n(A)$ and for all $n\in \N $,
\item[(4a)] There exists $c\ge 1$ such that \\
$\| (a_{lm})\| \le c\cdot \sup_{\alpha }\| (a_{lm}e_{\alpha })\| $
for all $(a_{lm})\in \M _n(A)$ and for all $n\in \N $,
\item[(5)] There exists $n\in \N $ such that
$\| (a_{lm})\| =\sup_{\alpha }\| (a_{lm}e_{\alpha })\| $ for all
$(a_{lm})\in \M _n(A)$,
\item[(5a)] There exist $c\ge 1$ and $n\in \N $ such that \\
$\| (a_{lm})\| \le c\cdot \sup_{\alpha }\| (a_{lm}e_{\alpha })\| $
for all
$(a_{lm})\in \M _n(A)$,
\item[(6)] $\| (a_{lm})\| =\sup \{ \| (a_{lm})(k_{lm})\| ;
(k_{lm})\in \M _n(K), \| (k_{lm})\| \le 1\} $ \\
for all $(a_{lm})\in \M _n(A)$ and for all $n\in \N $,
\item[(7)] There exists $n\in \N$ such that \\
$\| (a_{lm})\| =\sup \{ \| (a_{lm})(j_{lm})\| ;
(k_{lm})\in \M _n(K), \| (k_{lm})\| \le 1\} $ for all
$(a_{lm})\in \M _n(A)$,
\item[(8)] There exists $c\ge 1$ such that \\
$\| (a_{lm})\| \le c\cdot \sup \{ \| (a_{lm})(k_{lm})\| ;
(k_{lm})\in \M _n(K), \| (k_{lm})\| \le 1 \} $ \\
for all $(a_{lm})\in \M  _n(A)$ and for all $n\in \N $,
\item[(9)] There exist $c\ge 1$ and $n\in \N $ such that \\
$\| (a_{lm})\| \le c\cdot \sup \{ \| (a_{lm})(k_{lm})\| ;
(k_{lm})\in \M _n(K), \| (k_{lm})\| \le 1 \} $ for all
$(a_{lm})\in \M _n(A)$,
\item[(10')] A is a tight $\C -A$-essential extension of K,
\item[(11')] A is a $\C -A$-essential extension of K,
\item[(12')] A is a tight $A-A$-essential extension of K,
\item[(13')] A is an $A-A$-essential extension of K,
\item[(14)] $\M _{n}(K)$ is an essential ideal in $\M _{n}(A)$
for all $n\in \N $,
\item[(15)] $\M _{n}(K)$ is an essential ideal in $\M _{n}(A)$
for some $n\in \N $,
\item[(16)] $\M _n(A)$ is *-isomorphically embedded in
$\Ml (\M _n(K))$ for all $n\in \N $,
\item[(17)] $\M _n(A)$ is *-isomorphically embedded in
$\Ml (\M _n(K))$ for some $n\in \N $.
\end{itemize}
\end{corollary}

Thus, in the two-sided case, all the concepts of essentiality that we have 
introduced are equivalent.

From the above corollary and by observing ``$\M _n(J)$ is an essential left
ideal in $\M_n(A)$  $\Longleftrightarrow \M _n(JJ^*)$  is an essential
two-sided ideal in $\M _n(A)$''
(Definition\ref{essentialleftideals}), the following also holds.

\begin{corollary}
\label{JJ*}
Let A, J be as the assumptions in Theorem\ref{tight} and let $\{ u_{\beta }\} $
be a contractive approximate identity of JJ*. Then each statement of {\em (II)}
in Theorem\ref{tight} is equivalent to each of the following.
\begin{itemize}
\item[(1b)] I(JJ*) is completely isometrically isomorphic to I(A), \\
via a map that restricts to the identity on JJ*,
\item[(2b)] A is a tight $\C-\C$-essential extension of JJ*,
\item[(3b)] A is a tight $A-\C $-essential extension of JJ*,
\item[(3ab)] A is an $A-\C $-essential extension of JJ*,
\item[(4b)] $\| (a_{lm})\| =\sup_{\beta }\| (a_{lm}u_{\beta })\| $ for all
$(a_{lm})\in \M _n(A)$ and for all $n\in \N $,
\item[(4ab)] There exists $c\ge 1$ such that \\
$\| (a_{lm})\| \le c\cdot \sup_{\beta } \| (a_{lm}u_{\beta })\| $
for all $(a_{lm})\in \M _n(A)$ and for all $n\in \N $,
\item[(5b)] There exists $n\in \N $ such that
$\| (a_{lm})\| =\sup_{\beta }\| (a_{lm}u_{\beta })\| $ for all
$(a_{lm})\in \M _n(A)$,
\item[(5ab)] There exist $c\ge 1$ and $n\in \N $ such that
$\| a\| \le c\cdot \sup_{\beta }\| ae_{\beta }\| $ for all
$(a_{lm})\in \M _n(A)$,
\item[(6b)] $\| (a_{lm})\| =\sup \{ \| (a_{lm})(k_{lm})\| ;
(k_{lm})\in \M _n(JJ^*), \| (k_{lm})\| \le 1\} $ \\
for all $(a_{lm})\in \M  _n(A)$ and for all $n\in \N $,
\item[(7b)] There exists $n\in \N$ such that \\
$\| (a_{lm})\| =\sup \{ \| (a_{lm})(k_{lm})\| ;
(k_{lm})\in \M _n(JJ^*), \| (k_{lm})\| \le 1\} $ for all
$(a_{lm})\in \M_n(A)$,
\item[(8b)] There exists $c\ge 1$ such that \\
$\| (a_{lm})\| \le c\cdot \sup \{ \| (a_{lm})(k_{lm})\| ;
(k_{lm})\in \M _n(JJ^*), \| (k_{lm})\| \le 1 \} $ \\
for all $(a_{lm})\in \M  _n(A)$ and for all $n\in \N $,
\item[(9b)] There exist $c\ge 1$ and $n\in \N $ such that \\
$\| (a_{lm})\| \le c\cdot \sup \{ \| (a_{lm})(k_{lm})\| ;
(k_{lm})\in \M _n(JJ^*), \| (k_{lm})\| \le 1 \} $ for all
$(a_{lm})\in \M _n(A)$,
\item[(10'b)] A is a tight $\C -A$-essential extension of JJ*,
\item[(11'b)] A is a $\C -A$-essential extension of JJ*,
\item[(12'b)] A is a tight $A-A$-essential extension of JJ*,
\item[(13'b)] A is an $A-A$-essential extension of JJ*,
\item[(14b)] $\M _{n}(JJ^*)$ is an essential ideal in $\M _{n}(A)$
for all $n\in \N $,
\item[(15b)] $\M _{n}(JJ^*)$ is an essential ideal in $\M _{n}(A)$
for some $n\in \N $,
\end{itemize}
\end{corollary}

\section{Examples and applications.}
\label{Examples.}
It is easy to construct an example of a left but not two-sided ideal which
satisfies all the statements in Theorem\ref{tight}.

\begin{example}
{\em Let $A$ be any C*-algebra which properly contains an essential two-sided
ideal $K$. Set}
$$J:=\left(
\begin{array}{cc}
A & K \\
A & K
\end{array}
\right).$$
{\em Then $J$ is a left but not two-sided ideal in $\M _2(A)$ and
$\M _2(K)\subset J\subset \M _2(A)$, hence we can make
$I(\M _2(K))\subset I(J)\subset I(\M _2(A))$ with an injective envelope
$(I(\M _n(K)),\iota)$ of $\M _n(K).$ The identity map on $I(\M _2(K))$ extends
to a completely contractive linear map
$\phi : I(\M _2(A))\longrightarrow I(\M _2(K)).$ $\M _2(A)$ is a tight
$A-\C $-essential extension of $\M _2(K)$ by Corollary\ref{two-sided}, and
$I(\M _2(A))$ is a tight $A-\C $-essential extension of $\M _2(A)$ by
Theorem\ref{tightext}, so that $I(\M _2(A))$ is a tight $A-\C $-essential
extension of $\M _2(K)$. Hence by rigidity (Theorem\ref{tightext}), $\phi $ is 
a complete isometry since it fixes each element of $\iota (\M _2(K))$. Thus
$I(\M _2(A))=I(J)=I(\M _2(K))$.}
\end{example}

The following example shows that each statement in (II) in Theorem\ref{tight}
does not necessarily imply (5) in Theorem\ref{tight} or (5a) in
Theorem\ref{non-tight}. Especially, $J$ is an essential left ideal in $A$ is 
equivalent to saying that $A$ is a (tight) $\C -J$-essential extension of $J,$ 
but those do not necessarily imply that $A$ is a (tight) $A-\C $-essential extension of $J.$

\begin{example}
{\em Let $A:=\M _2$ and
$$J:=\left(
\begin{array}{cc}
\C & O \\
\C & O
\end{array}
\right).$$
Then $A$ is a C*-algebra and $J$ is a closed left ideal in $A$ with a
contractive right identity.
$e:=\left(
\begin{array}{cc}
1 & 0 \\
0 & 0
\end{array}
\right)$.
It is easy to see $JJ^*=A$, so especially, $JJ^*$ is an essential two-sided
ideal in $A$, hence $J$ is an essential left ideal in $A$ by
Definition\ref{essentialleftideals}. Let
$a:=\left(
\begin{array}{cc}
0 & 0 \\
0 & 1
\end{array}
\right),$ and let $(a_{lm})\in \M _{n}(A)$ be such that
$a_{11}=a,$ $a_{lm}=0$ if $l\ne 1$ or $m\ne 1.$ Then $\| (a_{lm})\| =1$, while 
$\| (a_{lm}e)\| =0.$}
\end{example}

\begin{remark}
{\em In general, if $A$ is a C*-algebra and $J$ is a closed left ideal in $A$,
and if $A$ is a (resp., tight) $A-\C $-essential extension of
$J$, and $J$ is properly contained in $A$, then $J$ does not have a (resp.,
contractive) right identity. In fact, suppose that $J$ has a (resp.,
contractive) right identity $e$. Then
$(a_{lm})\longmapsto (a_{lm}e) \; \forall (a_{lm})\in \M _n(A),
\forall n\in \N $ defines a completely bounded (resp., completely contractive)
left $A$-module map $\phi : A\longrightarrow J$ which is completely isometric
on $J$, so that $\phi $ is a completely bounded isomorphism (resp., a complete
isometry). But $0\ne (a-ae)\in $Ker$\phi $ for $a\in A\setminus J$, hence a
contradiction. Especially, for any finite dimensional left ideal 
$J\subsetneqq A,$ $A$ can not be a (tight) $A-\C $-essential extension of $J.$}
\end{remark}

We close this section with another application of a part of Theorem\ref{tight}.
As is defined in \cite{PSS}, for two C*-algebras $A$ and $B$ with 
$A,B \subset \B (\Hi)$ for some Hilbert space $\Hi,$ we say that $A$ 
{\em norms} $B$ when the following equation holds for each $n\ge 1$ and for 
each $X\in \M _n(B).$ 
$$\| X\| =\sup \{ \| RXC\| ; R\in Row_n(A), C\in Col_n(A), \|R\|,\|C\| \le 1\}
$$  Where $Row_n(A)$ and  $Col_n(A)$ are, respectively, row and column matrices
over $A.$ We show that, if $A$ is an essential two-sided ideal in $B,$ then $A$
norms $B.$ By Lemma 2.4 in \cite{PSS}, it suffices to show that, for each 
$n\ge 1$ and for each $X\in \M _n(B),$ 
$$\| X\| =\sup \{ \| XC\| ; C\in Col_n(A), \|C\| \le 1\}.$$
But $Col_n(A)$ is an essential left ideal in $\M _n(B),$ so the equation 
follows from Theorem\ref{tight} (14)$\implies$(6). 

\section{Conclusions and questions.}
(1) Compared with the tight case, the non-tight case is generally unknown. The 
difficulty in the non-tight case mainly comes from the lack of a rigidity 
result that does not need some module actions. Even in the Banach space case, 
there are a few deep results, but the entire picture is unclear. We do not know
the existence of a ``non-tight injective envelope'', namely a minimal injective
extension, for an arbitrary operator space. Also the lack of rigidity makes it 
difficult to connect essential extensions with even tight injective envelopes. 
As a result, we do not know  if either of the following implications is true:
$A$ is a tight $A-\C$-essential extension of $J$
$\stackrel{?}{\Longleftrightarrow }$ $A$ is an $A-\C$-essential extension of 
$J.$

(2) The implications 
(5)$\stackrel{?}{\implies }$(4)$\stackrel{?}{\implies }$(3) in 
Theorem\ref{tight} and
(5a)$\stackrel{?}{\implies }$(4a)$\stackrel{?}{\implies }$(3a) in 
Theorem\ref{non-tight} are still unknown.

(3) As one possible generalization of the results in this paper, one can 
consider replacing C*-algebras by operator algebras. In such a case,
the difficulty comes from the fact that we still do not know whether or not 
the
representation $\B (\Hi)$ of an operator $A-B$-bimodule
(\cite{CES}, \cite{B2}, \cite{BP}) is $A-B$-injective. Consequently, we do not 
know the
existence of any $A-B$-injective operator modules.
Such modules play an important role in
characterizing and constructing  essential extensions.

\vspace{5mm}

ACKNOWLEDGMENTS.
The first author wishes to express his gratitude to the second author, who is 
also his advisor, for introducing him to this subject and for constant 
encouragement. Both authors thank Professor David Peter Blecher for many 
useful helps.

\end{document}